\documentclass{article}

\usepackage[a4paper]{geometry}

\usepackage{amssymb}
\usepackage{graphicx}
\usepackage{mathptmx}

\newtheorem{theorem}{Theorem}

\newtheorem{corr}{Corollary}

\title{\textbf{LIBOR Interpolation And The HJM Model}\footnotemark\footnote{This research has been done at Vienna University of Technology and Dublin City University. The author gratefully acknowledges the Austrian Christian Doppler Society (CD-laboratory PRisMa), as well as Science Foundation Ireland (Edgeworth Center and FMC2) for their support. Special thanks go to Dr. Friedrich Hubalek for critical discussion and advice on the subject.}}

\author{\small{\textbf{Andreas Hula}}\\ \footnotesize{School Of Mathematics, Dublin City University,Dublin 9,Dublin,Ireland}\\ \footnotesize{(andreas.hula@dcu.ie)}}

\date{\today}

\begin{document}

\maketitle
\tableofcontents
\vspace{0.6 cm}
\textbf{HJM equations \and LIBOR market models \and L\'evy Processes \and Termstructure Interpolation } \\
\textbf{60G51, 91G30 , 60F99}
 \vspace{0.2 cm}
\begin{abstract}
We follow the lines of Musiela and Rutkowski \cite{MR97} and extend their interpolation method to models with jumps. Together with an extension method for the tenor structure of a given LIBOR market model (LMM) we get an infinite LIBOR termstructure.\\
Furthermore we present an argument why certain known exponential moment conditions on the HJM Model are necessary. The approach uses finite tenor LIBOR market models as approximation for the  HJM model, then extends and interpolates the tenor structure, relating it to the HJM structure.
\end{abstract}
\section{Introduction}
\label{intro}
This paper aims to extend a given finite forward LIBOR termstructure model to an infinite termstructure. Using the spot-LIBOR numeraire we enlarge a LIBOR market model beyond the given termstructure. Then we extend a technique for interpolation between LIBOR rates, first presented in \cite{MR97} to a situation were jumps are involved. Combining those methods to prove a summarizing theorem we get forward LIBOR dynamics for arbitrarily large maturities $T^*$ and for every maturity $T\in (0,T^*]$. 

\section{Building Extension And Interpolation}
\label{buildinter}
\subsection{Spot LIBOR rates}
\label{spot}
A spot measure by the definition of \cite{MR97} in discrete tenor LIBOR Market Model theory is given by the choice of numeraire $B(t,T_1)/B(0,T_1)$ for a given tenor structure.\\
We start from the dynamics under the proper forward-measure $\mathbb{P}_{T_{n+1}}$ for a finite discrete tenor model.(with tenor $\{ T_i | i=1,\dots,n,n+1  \}$)
\[
dL(t,T_n)=L(t_{-},T_n)(\lambda (t,T_n) c^{\frac12}_t dW^{n+1}_t +\int_{\mathbb{R}}(e^{\lambda ( t,T_n )x}-1 )(\mu - \nu_t^{n+1} )(dt,dx))
\]
and build the measure change as follows
\begin{equation}
(\frac{d\mathbb{P}_{T_1}}{d\mathbb{P}_{T_{n+1}}})_t=\frac{B(0,T_{n+1})B(t,T_1)}{B(0,T_{1})B(t,T_{n+1})}= \frac{F_B(t,T_1,T_{n+1})}{F_B(0,T_1,T_{n+1})}=\prod_{j=1}^n \frac{F_B(t,T_j,T_{j+1})}{F_B(t,T_j,T_{j+1})}
\end{equation}
yielding a Brownian Motion
\[
W_t^{1}=W_t^{n+1}-\int_0^t\sum_{j=1}^{n}\ell (s_{-},T_j)\lambda (s,T_j)c_s^{\frac12}ds \Rightarrow W_t^{n+1}=W_t^{1}+\int_0^t\sum_{j=1}^{n}\ell (s_{-},T_j)\lambda (s,T_j)c_s^{\frac12}ds
\]
and a compensator
\[
\nu^{1}=\prod_{j=1}^n \beta (t,x,T_j,T_{j+1}) \nu_t^{n+1} \Rightarrow \nu^{n+1}=\prod_{j=1}^n \frac1{\beta (t,x,T_j,T_{j+1})} \nu_t^1.
\]
The dynamics of $L(t,T_{n})$ under that measure $\mathbb{P}_{T_1}$ are then
\begin{equation}
dL(t,T_{n})=L(t_{-},T_n)\Big{(}\sum_{j=1}^n \ell(t_{-},T_j)\lambda (t,T_j)\lambda (t,T_n) c_t dt + 
\end{equation}
\[
 \int_{\mathbb{R}} (e^{\lambda (t,T_n)x}-1) (1-\prod_{j=1}^n \frac1{\beta (t,x,T_j,T_{j+1})}) \nu^1_t (dt,dx)+  \lambda (t,T_n)c_t^{\frac12} dW^1_t + \int_{\mathbb{R}} (e^{\lambda (t,T)x}-1) ( \mu-\nu^1_t )(dt,dx) \Big{)}.
\]

More generally starting from an arbitrary rate under its forward measure 
\[
dL(t,T_i)=L(t_{-},T_i)(\lambda (t,T_i) c^{\frac12}_t dW^{i+1}_s +\int_{\mathbb{R}}(e^{\lambda ( t,T_i )x}-1)(\mu - \nu_t^{i+1} )(dx,dt))
\]
The measure change becomes 
\begin{equation}
(\frac{d\mathbb{P}_{T_1}}{d\mathbb{P}_{T_{i+1}}})_t=\frac{B(t,T_1)}{B(t,T_{i+1})}=\frac{F_B(t,T_1,T_{i+1})}{F_B(0,T_1,T_{i+1})}=\prod_{j=1}^i\frac{F_B(t,T_j,T_{j+1})}{F_B(0,T_j,T_{j+1})}
\end{equation}
and the dynamics therefore
\begin{equation}
dL(t,T_i)=L(t_{-},T_i)\Big{(}\sum_{j=1}^i \ell (t_{-},T_j)\lambda (t, T_j)\lambda (t,T_i) c_t dt +
\end{equation}
\[
 \int_{\mathbb{R}}(e^{\lambda (t,T_i) x}-1) (1- \prod_{j=1}^i \frac1{\beta (t,x,T_j,T_{j+1})} ) \nu^1 (dx,dt) +c^{\frac12}_t \lambda (t,T_i)dW^1_t + \int_{\mathbb{R}} e^{\lambda (t,T_i)x}-1 ( \mu-\nu_t^1 )(dt,dx) \Big{)}
\]
by the same arguments concerning $W^1_t$ and $\nu^1_t$ as above.\\
Then therefore there exists an equivalent measure $\mathbb{P}_{T_1}$ to the other forward measures and the risk neutral measure, such that LIBOR-dynamics are of the form
\begin{equation}
dL(t,T_i)=L(t_{-},T_i)\Big{(}\sum_{j=1}^i \lambda (t,T_j)\ell (t_{-},T_j)\lambda (t,T_i) c_tdt +
\end{equation}
\[
 \int_{\mathbb{R}}(e^{\lambda (t,T_i)x}-1)(1-\prod_{j=1}^i\frac1{\beta (t,x,T_j,T_{j+1})}) \nu_t^{1}(dx,ds)+ \lambda (t,T_i)c_t^{\frac12} dW^1_t +\int_{\mathbb{R}} (e^{\lambda (t,T_i) x}-1) (\mu-\nu_t^1 ) (dx,dt) \Big{)}
\]
where
\begin{equation}
\nu_t^{i+1} = (\prod_{j=1}^i\frac1{(\ell (t_,T_j)(e^{\lambda (t,T_j)x}-1)+1)})\nu^1_t .
\end{equation}
So the LIBOR-rates in this case depend only on the rates modelled for shorter maturities.\\
Of course, a problem of this approach is, that $B(t,T_1)$ is essentially only defined on $[0,T_1]$ since afterwards the bond has matured.\\
Therefore we have to extend $B(t,T_1)$ beyond $T_1$ in such a way that it stays a semimartingale.\\
We will not go deeper into diverse possibilities for this, as we are primarily interested in an particular extension, which was introduced explicitly by Jamshidian.
\subsection{Spot-LIBOR Measure}
\label{spotLIBOR}
Introduced by Jamshidian in \cite{J1999}, was the numeraire
\begin{equation}
B^*(t)=\frac{B(t,T_{i(t)})}{B(0,T_1)}\prod_{j=1}^{i(t)-1}\frac{B(T_j,T_j)}{B(T_j,T_{j+1})}
\end{equation}
with $i(t)=\min \{ i:t\leq T_i  \}$.\\
The idea behind this is the following: We wish to extend a given $B(t,T_i)$ as explained above to time intervals beyond $[0,T_i]$. We do this by making $B(t,T_i)$ proportional to $B(t,T_{i+1})$ on $[T_i,T_{i+1}]$ , proportional to $B(t,T_{i+2})$ on $[T_{i+2},T_{i+3}]$ etc. That gives
\begin{equation}
B(t,T_i)=B(t,T_{j+1})\prod_{k=i}^j\frac{B(T_k,T_k)}{B(T_k,T_{k+1})},\qquad t\in [T_j,T_{j+1}),\qquad 1\leq i \leq j\leq n.
\end{equation}
For $B(t,T_1)$ this yields
\begin{equation}
B(t,T_1)=B(t,T_{i(t)})\prod_{j=1}^{i(t)-1}\frac{B(T_j,T_j)}{B(T_j,T_{j+1})}=B(t,T_{i(t)})\prod_{j=1}^{i(t)-1}(1+\delta_j L(T_j,T_j))\qquad \forall t \leq T_{n+1}.
\end{equation}
This can be interpreted as as the value of a bond from investing a given amount $B(0,T_1)$ at time $0$ at spot LIBOR rate $L(0,T_1)$ and at $T_1$ reinvesting the principal interest at the prevailing spot LIBOR rate $L(T_1,T_1)$ and so on.\\
The numeraire $B^*(t)$ is then given through the extended $B(t,T_1)$ as
\[
B^*(t)=\frac{B(t,T_1)}{B(0,T_1)}
\]
We will show the following
\begin{theorem}[Spot-LIBOR Dynamics]
There is a measure denoted by $\mathbb{P}_{Ls}$ given through the numeraire $B^*(t)$, equivalent to the forward-measures and the risk neutral measure such that the dynamics of the LIBOR rates for a given tenor structure are 
\begin{equation}
L(t,T_s)=L(t_{-},T_s)(\sum_{j=i(t)}^s \ell (t_{-},T_j)\lambda (t,T_j)\lambda (t, T_s) c_t dt +\int_{\mathbb{R}}(e^{\lambda (t,T_s) x}-1)
\end{equation}
\[
 (1- \prod_{j=i(t)}^s \frac1{\beta (t,x,T_j,T_{j+1})} ) \nu^{i(t)}_t (dx,dt) + \lambda (t,T_s)c_t^{\frac12} dW^{i(t)}_t + \int_{\mathbb{R}} (e^{\lambda (t,T_s)x}-1) ( \mu - \nu^{i(t)}_t )(dt,dx) ).
\]
\end{theorem}
$Proof:$\\
To understand the LIBOR-market-model under the measure induced by that numeraire, we first look at the first time interval under consideration $[0,T_1]$ and see our spot-measure as discussed above. How to continue for the other time-intervals? To answer this, we look at the form of the measure change for each time interval 
\begin{equation}
\frac{d\mathbb{P}_{Ls}}{d\mathbb{P}_{T_{s+1}}}=\frac{B(t,T_{i(t)-1})}{B(0,T_1)}\prod_{j=1}^{i(t)-1}\frac{B(T_j,T_j)}{B(T_j,T_{j+1})}\frac{B(0,T_{s+1})}{B(t,T_{s+1})}=
\end{equation}
\[ 
\frac{B(t,T_{i(t)})}{B(t,T_{s+1})}\prod_{j=0}^{i(t)-1}\frac{B(T_j,T_j)}{B(T_j,T_{j+1})}= 
\frac{F_B(t,T_{i(t)},T_{s+1})}{F_B(0,T_{i(t)},T_{s+1})}\prod_{j=1}^{i(t)-1}F_B(T_j,T_j,T_{j+1})=\frac{F_B(t,T_{i(t)},T_{s+1})}{F_B(0,T_{i(t)},T_{s+1})}C
\]
So for each interval $t\in (T_i,T_{i+1}]$ our numeraire is then $B(t,T_{i+1}) C$. We therefore can express the Spot-LIBOR numeraire dynamics by a sequence of forward-measure dynamics.\\
We compute measure changes accordingly
\begin{equation}
(\frac{d\mathbb{P}_{Ls}}{d\mathbb{P}_{T_{s+1}}})_t=\prod_{j=i(t)}^{s}F_B(t,T_{j},T_{j+1})F_B(0,T_{j},T_{j+1})
\end{equation}
Inserting the resulting equation for the Brownian Motion
\begin{equation}
W_t^{s+1}=W_t^{i(t)}+\int_0^t \sum_{j=i(t)}^s \lambda (t,T_j)\ell (u_{-},T_j)c_u^{\frac12}du
\end{equation}
and the compensator
\begin{equation}
\nu_t^{s+1}=\prod_{j=i(t)}^s \frac1{\beta (t,x,T_j,T_{j+1})}\nu_t^{i(t)}
\end{equation}
into the forward dynamics
\begin{equation}
dL(t,T_s)=L(t_{-},T_s)\Big{(}\lambda (t,T_s)c_t^{\frac12}dW_t^{s+1}+\int_{\mathbb{R}}(e^{\lambda (t,T_s)x}-1)(\mu-\nu_t^{s+1} ) (dt,dx) \Big{)}
\end{equation}
This yields then the following dynamics
\begin{equation}
L(t,T_s)=L(t_{-},T_s)(\sum_{j=i(t)}^s \ell (t_{-},T_j)\lambda (t,T_j)\lambda (t,T_s) c_t dt +\int_{\mathbb{R}}(e^{\lambda (t,T_s) x}-1)
\end{equation}
\[
 (1- \prod_{j=i(t)}^s \frac1{\beta (t,x,T_j,T_{j+1})} ) \nu^{i(t)}_t (dx,dt) + \lambda (t,T_s)c_t^{\frac12} dW^{i(t)}_t +
 \]
 \[
 \int_{\mathbb{R}} (e^{\lambda (t,T_s)x}-1) ( \mu - \nu^{i(t)}_t )(dx,dt) )\qquad \forall t \in [T_{i(t)-1},T_{i(t)}]
\]
$\quad\Box$
Those dynamics, where each rate is dependent only on finitely many( already calculated) rates, form the basis for our extension of a given model to an infinite time horizon.\\
This is especially interesting since in an HJM framework we would assume models to be defined for arbitrary large maturities.\\
For the time being we can only work on a discrete time-grid, but this problem can be solved by "filling" the gaps as we show in the section on continuous tenors.
\subsection{Extending The Tenor}
\label{Tenorext}
\begin{theorem}[LIBOR-Extension]\label{ext}
For any given finite tenor-structure $\{ T_i \}_{i=1}^{n+1}$,  strictly decreasing, positive initial term structure $(B(0,T_i))_{i=1}^{n+1}$ and volatility functions $\{\lambda (.,T_i) \}_{i=1}^{n}$ and a corresponding LIBOR-Market-Model $(L(.,T_i))_{i=1}^n$ we may choose positive functions $\{\lambda (.,T_i)\}_{i=n+1}^{\infty}$ such that the integrability condition is fulfilled and from that obtain a unique extension of our model $\{ L(.,T_i) \}_{i=1}^{\infty}$ by demanding each of our LIBOR-rate processes fulfills the finite dimensional SDE
 \begin{equation}
dL(t,T_s)=L(t_{-},T_{s})( \sum_{j=i(t)}^{s}\ell (t_{-},T_j)\lambda (t,T_j)\lambda (t,T_s)c_t dt+\lambda (t,T_s)c_t^{\frac12}dW_t^{i(t)}+
\end{equation}
\[
 \int_{\mathbb{R}}(e^{\lambda (t,T_{s})x}-1)(1-\prod_{j=i(t)}^{s} \frac1{ \beta (t,x,T_j,T_{j+1})}) \nu^{i(t)}(dx,dt) +\int_{\mathbb{R}}(e^{\lambda (t,T_{s})x} -1 )(\mu-\nu_t^{i(t)} )(dt,dx))
\] 
given an initial condition $L(0,T_s)$.
\end{theorem}
$Proof$:\\
We assume we are given a finite tenor-structure and we are working under the Spot-LIBOR measure.\\
Say we add another point to the tenor structure $T_{n+2} > T_{n+1}$ and $T_{n+2}-T_{n+1}=\delta_{n+1}$. We have the following relation between Brownian Motion under $\mathbb{P}_{T_{n+2}}$ for an arbitrary forward measure and Brownian Motion for $\mathbb{P}_{Ls}$ :
\begin{equation}
W_t^{n+2}=W_t^{i(t)}+\int_0^t \sum_{j=i(s)}^{n+1}\frac{\delta_j L(u_{-},T_j)}{1+\delta_j L(u_{-},T_j)}\lambda (u,T_j)c_u^{\frac12}du\qquad u_{-} \in [T_{i(s)-1},T_i(s)) 
\end{equation}
and the compensator
\begin{equation}
 \nu^{n+2}=\prod_{j=i(t)}^{n+1}\frac1{\beta (t,x,T_j,T_{j+1})} \nu_t^{i(t)}\qquad \forall t\in [T_{i-1},T_i].
\end{equation}
We can therefore write down a SDE for $L(t,T_{n+1})$ under $\mathbb{P}_{T_{n+2}}$
\begin{equation}
dL(t,T_{n+1})=L(t_{-},T_{n+1})( \sum_{j=i(t)}^{n+1}\ell (t_{-},T_j)\lambda (t,T_j)\lambda (t,T_{n+1})c_t dt+ c_t^{\frac12}dW^i(t)_t+
\end{equation}
\[
 \int_{\mathbb{R}}(e^{\lambda (t,T_{n+1})x}-1)(1-\prod_{j=i(t)}^{n+1} \frac1{ \beta (t,x,T_j,T_{j+1})}) \nu^{i(t)}(dx,dt) + \int_{\mathbb{R}}(e^{\lambda (t,T_{n+1})x} -1 )(\mu-\nu_t^{i(t)} )(dx,dt))
\]
for a a priori unspecified positive bounded function $\lambda (t,T_{n+1})$( to be determined through calibration for instance).$\quad\Box$
Since we know how to switch between forward and spot-LIBOR measures, it does not matter under which measure we originally specify our LIBOR-Market-Model. \\

Obviously we can repeat this procedure, choosing a positive function $\lambda (t,T_{n+2})$, a new point in time $T_{n+3}> T_{n+2}$ getting a well defined, solvable( finite-dimensional, with Lipschitz-Coefficients if the $\lambda (.,T_i)$ are chosen that way) SDE for $L(t,T_{n+2})$.\\
Therefore, if we extend our tenor-structure, to an arbitrarily large (even countably infinite) set of time points $\{T_i \}_{i=1}^{\infty}$, we get for any possible rate
\begin{equation}
dL(t,T_s)=L(t_{-},T_{s})( \sum_{j=i(t)}^{s}\ell (t_{-},T_j)\lambda (t,T_j)\lambda (t,T_s)c_s dt+\lambda(t,T_s) c_t^{\frac12}dW_t^{i(t)}+
\end{equation}
\[
 \int_{\mathbb{R}}(e^{\lambda (t,T_{s})x}-1)(1-\prod_{j=i(t)}^{s} \frac1{ \beta (t,x,T_j,T_{j+1})}) \nu^{i(t)}(dx,dt) +\int_{\mathbb{R}}(e^{\lambda (t,T_{s})x} -1 )(\mu-\nu_t^{i(t)} )(dt,dx))
\] 
which is a finite dimensional SDE for any fixed $T$, dependent only on already calculated rates and therefore solvable with purely finite-dimensional methods. 
Such an infinite discrete tenor structure may serve as a skeleton for a full continuous tenor term structure model. We would have to fill the gaps between the tenor points. We will address this question in the section on continuous tenors.
We wish to extend the construction of Musiela and Rutkowski \cite{MR97} to semimartingale driven LIBOR-Models so as to get a "full" tenor-structure in that case as well.
\subsection{Construction Concept}
\label{constr}
In analogy to the work of Musiela and Rutkowski \cite{MR97}, we wish to "fill the gaps" between the discrete tenor dates $\{T_i\}_{i\in 1,\dots, n+1}$. For that we assume an equidistant tenor-time-grid. \\
At first we will assume to be working up to a terminal maturity $T_{n+1}$ and wish to specify the dynamics of $L(t,T)$ for all $T\in [ 0,T_{n+1} ]$. Later on arbitrary maturities $T\in\mathbb R_+$ will be considered. \\
As in \cite{MR97} we use backward induction for this
\begin{enumerate}
\item
First, we define a forward LIBOR-market model on a given equidistant discrete grid $T_i=i\delta$.
\item
Secondly, numeraires for the interval $( T_n,T_{n+1} )$. We have values for the spot-LIBOR numeraire at $T_n$ and $T_{n+1}$, in short $B^*(T_{n})$ and $B^*(T_{n+1})$. Both $B^*(T_n)$ and $B^*(T_{n+1})$ are $\mathcal{F}_{T_{n}}$ measurable random variables.\\
We define a spot martingale measure through $\frac{d\mathbb{P}_{Ls}}{d\mathbb{P}_{T_{n+1}}}=B^*(T_{n+1})B(0,T_{n+1})$.\\
We attempt to satisfy intial conditions in our model for the interpolated rates via
a function $\gamma: [ T_n,T_{n+1} ]\rightarrow [0,1]$ such that $\gamma(T_n)=0$ and $\gamma(T_{n+1})=1$ and the process
\[
\log B^*(T)=(1-\gamma (T) )\log B^*(T_n) +\gamma (T) \log B^*(T_{n+1}),\qquad \forall T \in [T_n,T_{n+1}],
\]
satisfies $B(0,t)=\mathbb{E}_{\mathbb{P}_{Ls}}(1/B_t^*)$ for every $T\in[T_n,T_{n+1}]$. We have that $0< B^*(T_n)<B^*(T_{n+1})$ and $B(0,t),t\in [T_n,T_{n+1}]$ is assumed to be a strictly decreasing function, so such a $\gamma$ exists and is unique.
\item
Thirdly, given the spot-LIBOR numeraires $B^*(t)$ for all $t\in [T_n,T_{n+1}]$ the forward measure for any date $T\in (T_n,T_{n+1})$ can be defined by the formula 
\[
\frac{d\mathbb{P}_T}{d\mathbb{P}_{Ls}}=\frac{1}{B^*(T) B(0,T)}.
\]
If we use this and the definition of our spot martingale measure, we get
\[
\frac{d\mathbb{P}_{T}}{d\mathbb{P}}=\frac{d\mathbb{P}_T}{d\mathbb{P}_{Ls}}\frac{d\mathbb{P}_{Ls}}{d\mathbb{P}}=\frac{B^*(T_{n+1})B(0,T_{n+1})}{B^*(T) B(0,T)}
\]
which gives for every $T\in [T_n,T_{n+1}]$
\[
\frac{d\mathbb{P}_T}{d\mathbb{P}} |_{\mathcal{F}_t} =\mathbb{E}_{\mathbb{P}}(\frac{B^*(T_{n+1}) B(0,T_{n+1})}{B^*(T) B(0,T)}| \mathcal{F}_t)
\]
Using stochastic exponentials to describe this we get
\[
\frac{d\mathbb{P}_T}{d\mathbb{P}} |_{\mathcal{F}_t}=\frac{B(0,T_{n+1})}{B(0,T)}\mathcal{E}_t\Big{(}-\int_0^{.} \alpha (u,T,T_{n+1})c_u^{\frac12}dW_u^{n+1}+\int_0^{.}\int_{\mathbb{R}}(\beta(u,x,T,T_{n+1}) -1)(\mu-\nu_t^{n+1} )(dx,du) \Big{)}
\]
which we use to describe the forward volatility $\alpha (t,T,T_{n+1})$ for any maturity $T\in (T_n,T_{n+1})$. We get a $\mathbb{P}_T$ Wiener process $W^T$ and a $\mathbb{P}_T$ compensator for the jump-part. Given those ingredients we define the forward LIBOR rate process $L(t,T-\delta)$ for arbitrary $T\in (T_n,T_{n+1})$ by setting 
\[
dL(t,T-\delta )=L( t_{-},T-\delta )\Big{(}\lambda (t,T-\delta )c_t^{\frac12}dW_t^T+\int_{\mathbb{R}} (e^{\lambda (t,T-\delta )x}-1) (\mu-\nu_t^T )(dt,dx)\Big{)}
\]
with usual initial condition
\[
L(0,T-\delta )=\delta^{-1}(\frac{B(0,T -\delta )}{B(0,T)}-1).
\]
Finally we know
\[
\alpha (t,T_n,T_{n+1})=\ell (t_{-},T) \lambda (t,T_n)
\]
and
\[
\beta (t,x,T_n,T_{n+1})=\ell (t_{-},T)( e^{(\lambda (t,T_n) x)}-1)+1
\]
and thus we are able to define the forward measure for the date $T$.\\
To define forward probability measures $\mathbb{P}_{U}$ and the corresponding driving processes for all maturities $U\in (T_{n-1},T_n)$ we put
\[
\alpha (t,U,T)=\alpha(t, T-\delta ,T)= \frac{\delta L(t,T-\delta}{1+\delta L(t,T-\delta )} \lambda(t, T-\delta)
\]
and
\[
\beta (t,x,U,T)=\beta (t,x,T -\delta ,T)=\ell (t_{-},T) (e^{(\lambda (t,T -\delta )x)}-1)+1
\]
with $U=T - \delta $ such that $T=U+\delta$ belongs to $(T_{n},T_{n+1})$.\\
The relations between those coefficients are derived from the necessary relations between forward measure changes( see the section on forward modeling).\\
The coefficient $\alpha(t,U,T_{n+1})$ is calculated through
\[
\alpha (t,U,T_{n+1})=\alpha(t,U,T)-\alpha (t,T,T_{n+1}),\qquad \forall t\in [0,T -\delta ].
\]
For the jump part, $\beta(t,x,U,T_{n+1})$ is calculated through
\[
\beta(t,x,U,T_{n+1})=\beta(t,x,U,T) \beta(t,x,T,T_{n+1}),\qquad \forall t\in [0,T- \delta ].
\]
\end{enumerate}
Continuing this Backward construction, we get a continuous tenor LIBOR model.\\
Since we construct a family of forward measures, we can construct a family of forward processes $F(t,T_{n+1},T)$ which fulfill the SDE
\[
dF(t,T,T_{n+1})=F(t_{-},T,T_{n+1})(\alpha (t,T,T_{n+1})c_t^{\frac12}dW^{n+1}_t+\int_{\mathbb{R}}(\beta(t,x,T,T_{n+1})-1)(\mu_t-\nu^{n+1}_t)(dx,dt)).
\]
By construction we have that from those forward processes we get a family of bond prices $B(t,T)$ by $B(t,T):=F(t,T,t)$. The family of bond prices obtained thus always satisfies the weak no-arbitrage condition.\\
Now we transform the equations above to the spot-LIBOR measure.
\subsection{Spot-LIBOR Interpolation}
\label{spotLIBinter}
Now we want to carry out interpolation for a model given under the spot-LIBOR measure.\\
We assume a finite equidistant tenor-structure and 
 the spot-LIBOR dynamics
\begin{equation}
dL(t,T_s)=L(t_{-},T_s)\Big{(}\lambda (t,T_s)\sum_{j=i(t)}^s \lambda (t,T_j)\ell (t_{-},T_j)c_tdt+ \lambda (t,T_s)c_t^{\frac12}dW_t^{i(t)}+
\end{equation}
\[
 \int_{\mathbb{R}}(e^{\lambda (t,T)x}-1)(1-\prod_{j=i(t)}^s\frac1{\beta (t,x,T_j)})\nu_t^{i(t)}(dx,dt) +\int_{\mathbb{R}}(e^{\lambda (t,T)x}-1)(\mu-\nu_t^{i(t)})(dx,dt)  \Big{)}
\]
for all $T_s\in \{ T_i \}_{i=1}^{n+1}$.\\
We start in the interval $(T_{n},T_{n+1})$. We define the spot-LIBOR numeraire process $B(T)$ for all $T \in (T_{n},T_{n+1})$ just as in the section on forward interpolation above and assume a positive bounded function $\lambda (t,T_s)$.\\
We can calculate a change of measure from $\mathbb{P}_{T_{n+1}}$ to $\mathbb{P}_{T+\delta}$ as in the section above as well
\begin{equation}
\frac{d\mathbb{P}_{T+\delta}}{d\mathbb{P}_{T_{n+1}}}=\mathcal{E}(\int_0^t\alpha (s,T+\delta ,T_{n+1})c_t^{\frac12}dW_s^{n+1}+\int_0^t\int_{\mathbb{R}}(\beta (s,x,T+\delta ,T_{n+1})-1)(\mu-\nu_s^{n+1})(ds,dx))
\end{equation}
We know therefore, that a Brownian Motion for the forward measure $\mathbb{P}_{T+\delta}$ is given as
\[
W^{T+\delta}_t=W_t^{n+1}+\int_0^t \lambda (s,T)\ell (s,T)c^{\frac12}_s ds\qquad \forall t \in [0,T_{n+1}]
\]
in terms of $\mathbb{P}_{T_{n+1}}$ with the compensator being
\[
\nu_t^{T+\delta}=\nu_t^{n+1}\frac1{\beta (t,x,T_{n+1},T+\delta)}=\nu_t^{n+1}\frac1{\ell (t_{-},T)(e^{\lambda (t,T)x}-1)+1}.
\]
A forward LIBOR-rate for $T\in (T_n,T_{n+1}]$ has to fulfill
\[
dL(t,T)=L(t_{-},T)\Big{(} \lambda (t,T)c_t^{\frac12} dW_t^{T+\delta} +\int_{\mathbb{R}}(e^{\lambda (t,T)x}-1)(\mu-\nu^{T+\delta}_t)(dt,dx) \Big{)}
\]
under its proper forward measure $\mathbb{P}_{T+\delta}$.\\
Under $\mathbb{P}_{T_{n+1}}$ we then get $W_t^{T+\delta}=W_t^{n+1}+\int_0^t \lambda (t,T)\ell (s_{-},T)c_s^{\frac12}ds$ and $\nu_t^T=\nu_t^{n+1}\frac1{\beta(t,x,T_{n+1},T+\delta)}$.
From this we get
\begin{equation}
dL(t,T)=L(t_{-},T)\Big{(}\lambda (t,T)c_t^{\frac12}dW_t^{n+1}+\lambda (t,T)^2 c_t\ell (t_{-},T_{n+1})dt+
\end{equation}
\[
 \int_{\mathbb{R}}(e^{\lambda (t,T)x}-1)(1-\frac1{\beta (t,x,T_{n+1},T+\delta)})\nu_t^{n+1}+\int_{\mathbb{R}}(e^{\lambda (t,T)x}-1) (\mu-\nu_t^{n+1}) \Big{)}.
\]
in terms of $\mathbb{P}_{T_{n+1}}$.\\
In terms of the spot-LIBOR measure we get
\begin{equation}
dL(t,T)=L(t_{-},T)\Big{(}\sum_{j=i(t)}^{i(T)-1}\lambda(t,T_j)\lambda (t,T) c_t\ell (t_{-},T_j)dt+\ell(t_{-},T)c_t\lambda (t,T)^2)dt+
\end{equation}
\[
 \int_{\mathbb{R}}(e^{\lambda (t,T)x}-1)(1-\prod_{j=i(t)}^{i(T)-1} 
\frac1{\beta (t,x,T_j,T_{j+1})}) \nu^{i(t)}_t (dt,dx)+\int_{\mathbb{R}}(e^{\lambda (t,T)x}-1)(1-\frac1{\beta (t,x,T_{n+1},T+\delta )})\nu_t^{n+1}(dt,dx) +
\]
\[ 
\lambda (t,T)c_t^{\frac12}dW_t^{n+1}+\int_{\mathbb{R}}(e^{\lambda (t,T)x}-1)(\mu-\nu_t^{n+1})(dt,dx)  \Big{)}.
\]
So the interpolation on the interval $T\in (T_n,T_{n+1}]$ is well defined and we have that
\begin{equation}
 W_t^T=W_t^{n+1}+\int_0^t \lambda (s,T)\ell (s_{-},T)c_s^{\frac12}ds =W_t^{i(T)}+\int_0^t \lambda (s,T)\ell (s_{-},T)c_s^{\frac12}ds
\end{equation}
and
\begin{equation}
\nu_t^{T}=\nu_t^{i(T)}\frac1{\beta (t,x,T_{i(t)},T + \delta)}.
\end{equation}
Once we have the interpolated model for a whole interval, we use the relation 
\[
\frac{d\mathbb{P}_{T -k\delta}}{d\mathbb{P}_{T-(k-1)\delta}}=\frac{F_B(t,T-k\delta ,T-(k-1)\delta)}{F_B(0,T-k\delta ,T-(k-1)\delta )}\qquad\forall k \leq i(T)
\]
 to get Brownian Motions
\[
W_t^{T-(k-1)\delta}=W_t^{i(T-k\delta )}+\int_0^t \lambda (t,T-k\delta )\ell (s,T-k\delta) c_t^{\frac12} ds\qquad \forall -\infty <k \leq i(T)-1
\]
 and Compensators
 \[
\nu_t^{T-k\delta}=\nu_t^{i(T-k\delta)}\frac1{\beta (t,x,i(T-k\delta) , T-(k-1)\delta )}\qquad \forall -\infty <k \leq i(T)-1
 \]
 for all remaining maturities $T-k\delta\in[T_n-k\delta,T_{n+1}-k\delta]$ and our interpolated processes become solutions of
 \begin{equation}
 dL(t,T)=L(t_{-},T)\Big{(}\sum_{j=i(t)}^{i(T)-1}\lambda (t,T)\lambda(t,T_j)c_t \ell (t_{-},T_j)dt+\ell(t_{-},T)c_t\lambda (t,T)^2)dt+
 \end{equation}
 \[
 \int_{\mathbb{R}}(e^{\lambda (t,T)x}-1)(1-\prod_{j=i(t)}^{i(T)-1} \frac1{\beta (t,x,T_j,T_{j+1})}) \nu^i(t)_t(dt,dx) 
 +
\]
\[ 
\int_{\mathbb{R}}(e^{\lambda (t,T)x}-1)(1-\frac1{\beta (t,x,i(T)-1,T+\delta)})\nu_t^{i(T)}(dt,dx) + 
\lambda (t,T)c_t^{\frac12}dW_t^{n+1}+\int_{\mathbb{R}}(e^{\lambda (t,T)x}-1)(\mu-\nu_t^{i(T)})(dt,dx)  \Big{)}.
 \]
For an arbitrary starting interval the procedure works as follows:
\begin{enumerate}
\item
We look at $T\in (T_k,T_{k+1})$. We want to define the dynamics of $L(t,T)$ in an arbitrage-free way for all $T\in (T_k,T_{k+1})$. For that, we interpolate between two spot-LIBOR numeraires.
\[
\log B(T)^*=(1-\gamma(t) )\log B^*_{T_{k}} +\gamma(t) \log B^*_{T_{k+1}},\qquad \forall T \in [T_k,T_{k+1}],
\]
\item
We determine the measure change between $\mathbb{P}_{T_{k+1}}$and $\mathbb{P}_{T+\delta}$.
\begin{equation}
\frac{d\mathbb{P}_{T+\delta}}{\mathbb{P}_{T_{k+1}}}=\mathcal{E}(\int_0^t\alpha (s,T,T_{k+1})c_t^{\frac12}dW_s^{k+1}+\int_0^t\int_{\mathbb{R}}(e^{\lambda (s,T)}-1)(\mu-\nu_s^{k+1})(ds,dx))
\end{equation}
\item
We determine the Brownian Motion and the compensator for $\mathbb{P}_{T+\delta}$ in terms of the forward measure $\mathbb{P}_{T_{k+1}}$:
\begin{equation}
 W_t^T=W_t^{k+1}+\int_0^t \lambda (s,T)\ell (s_{-},T)c_s^{\frac12}ds
\end{equation}
and
\begin{equation}
\nu_t^{T}=\nu_t^{k+1}\frac1{\beta (t,x,T_{k+1},T+\delta)}.
\end{equation}
\item
From this we can determine the dynamics of $L(t,T)$ under the spot-LIBOR measure:
 \begin{equation}
 dL(t,T)=L(t_{-},T)\Big{(}\lambda (t,T)(\sum_{j=i(t)}^{k}\lambda (t,T_j)\ell (t_{-},T_j)c_t+\ell (t_{-},T)c_t\lambda (t,T))dt+
 \end{equation}
 \[
 \int_{\mathbb{R}}(e^{\lambda (t,T)x}-1)(1-\prod_{j=i(t)}^{k+1} 
\frac1{\beta (t,x,T_j,T_{j+1})}) \nu^{i(t)}_t(dt,dx) 
 +\int_{\mathbb{R}}(e^{\lambda (t,T)x}-1)(1-\frac1{\beta (t,x,T_{k+1},T+\delta)})\nu_t^{k+1}(dt,dx) +
\]
\[ 
\lambda (t,T)c_t^{\frac12}dW_t^{i(t)}+\int_{\mathbb{R}}(e^{\lambda (t,T)x}-1)(\mu-\nu_t^{i(t)})(dt,dx)  \Big{)}.
 \]
\end{enumerate}
From one fully determined interval we can determine the LIBOR-rate process dynamics of any other by consequence of the forward measure changes to be 
 \begin{equation}
 dL(t,T)=L(t_{-},T)\Big{(}\lambda (t,T)(\sum_{j=i(t)}^{i(T)-1}\lambda (t,T_j)c_t\ell (t_{-},T_j)+\ell (t_{-},T)\lambda (t,T)c_t)dt+ 
\end{equation}
\[ 
 \int_{\mathbb{R}}(e^{\lambda (t,T)x}-1)(1-\prod_{j=i(t)}^{i(T)-1} 
\frac1{\beta (t,x,T_j,T_{j+1})}) \nu^{i(t)}_t(dt,dx) 
 +
\]
\[
 \int_{\mathbb{R}}(e^{\lambda (t,T)x}-1)(1-\frac1{\beta (t,x,i(T),T+\delta)})\nu_t^{i(T)-1}(dt,dx) + 
\lambda (t,T)c_t^{\frac12}dW_t^{i(t)}+\int_{\mathbb{R}}(e^{\lambda (t,T)x}-1)(\mu-\nu_t^{i(t)})(dt,dx)  \Big{)}.
 \]
Hereby we have also shown, that our method does not depend on the particular choice of the starting interval( since every interval yields the same SDE). Through induction the restriction to a finite tenor structure is not necessary.
\section{Existence Of LIBOR-Term Structure Models}
\label{LIBterm}
We gather the results we have derived so far in the following theorem
\begin{theorem}\label{LIBORTERM}
Given a equidistant discrete tenor structure $\{ T_i \}_{i\in I}$( $I$ possibly infinite), volatility functions $\{ \lambda (t,T_i)\}_{i\in \mathbb{N}}$, an initial strictly positive, strictly decreasing term-structure $\Big{(}B(0,T)\Big{)}$ and a driving process
\[
X_t:=\int_0^t b(s,T_1)ds +\int_0^t c_s^{\frac12}dW_s^1+\int_0^t\int_{\mathbb{R}}x (\mu-\nu_s^1)(ds,dx)
\]
fulfilling
\begin{equation}\label{cond1}
\int_0^t (|b(s,T_1)|+c_s )ds<\infty \qquad \forall t \in\mathbb{R}_+
\end{equation}
as well as 
\begin{equation}\label{cond2}
\int_0^{\infty} \int_{|x|\geq 1} \exp (ux)F_s(dx) ds< \infty\qquad u < M, M\geq \sum_{i\in\mathbb{N}}|\lambda (.,T_i)|,M<\infty
\end{equation}
and
\begin{equation}\label{cond3}
\int_0^{\infty}\int_{\mathbb{R}} (x^2 \wedge 1)F_s(dx)ds<\infty
\end{equation}
then there is a LIBOR termstructure $\{ L(t,T) \}_{t\leq T, T\in\mathbb{R}_+}$ fulfilling
\begin{equation}
 dL(t,T)=L(t_{-},T)\Big{(}\lambda (t,T)(\sum_{j=i(t)}^{i(T)-1}\lambda(t,T_j) c_t\ell (t_{-},T_j)dt+\ell (t_{-},T)c_t\lambda (t,T))dt+
\end{equation}
\[ 
 \int_{\mathbb{R}}(e^{\lambda (t,T)x}-1)(1-\prod_{j=i(t)}^{i(T)-1} 
\frac1{\beta (t,x,T_j,T_{j+1})}) \nu^{i(t)}_t(dt,dx) 
 +\int_{\mathbb{R}}(e^{\lambda (t,T)x}-1)(1-\frac1{\beta (t,x,T_{i(T)},T+\delta)})\nu_t^{i(T)}(dt,dx) +
\]
\[ 
\lambda (t,T)c_t^{\frac12}dW_t^{i(t)}+\int_{\mathbb{R}}(e^{\lambda (t,T)x}-1)(\mu-\nu_t^{i(t)})(dt,dx)  \Big{)}.
\]
 for all $T\in\mathbb{R}_+$.
\end{theorem}
$Proof$:\\
\begin{enumerate}
\item
Start with a finite tenor structure $\{ T_i \}_{i=1}^{n+1}$.
\item
Define a LIBOR-Market-Model under the Spot-LIBOR measure $\mathbb{P}_{Ls}$ as solution to the corresponding SDE's.
\begin{equation}
dL(t,T_s)=L(t_{-},T_s)\Big{(} \lambda (t,T_s)\sum_{j=i(t)}^s\lambda (t,T_j)\ell (t_{-},T_j)c_tdt +\lambda (t,T_s )c_t^{\frac12}dW_t^{i(t)}+
\end{equation}
\[
 \int_{\mathbb{R}}(e^{\lambda (t,T_s)x}-1)(1-\prod_{j=i(t)}^s\frac1{\beta (t,x,T_j,T_{j+1})} )\nu_t^{i(t)}(dt,dx)+\int_{\mathbb{R}}(e^{\lambda (t,T_s)x}-1)(\mu-\nu_t^{i(t)})(dt,dx) \Big{)}.
\]
 Since there are only finitely many factors entering into each equation, we have the usual existence and uniqueness theorems in finite dimension.
\item 
Interpolate between every two tenor points. This is possible through arbitrage free interpolation under the spot-LIBOR measure. See section ~\ref{spotLIBinter}
\item
Extend this LIBOR-Market-Model to an infinite tenor structure as in the LIBOR extension theorem. This 
is well defined, as shown in the theorem ~\ref{ext}
\item
Interpolate/extend for every extension for the tenor-grid, the interpolated LIBOR-rate dynamics in between. 
\item
Interpret the resulting family of LIBOR-rate processes $\{L(t,T)\}_{t,T}$ as a solution of an infinite dimensional problem
\begin{equation}
 dL(t,T)=L(t_{-},T)\Big{(}\sum_{j=i(t)}^{i(T)-1}\lambda(t,T_j)c_t\ell (t_{-},T_j)dt+\ell(t_{-},T)c_t\lambda (t,T)^2)dt+
 \end{equation}
 \[
 \int_{\mathbb{R}}(e^{\lambda (t,T)x}-1)(1-\prod_{j=i(t)}^{i(T)-1} 
\frac1{\beta (t,x,T_j,T_{j+1})}) \nu^{i(t)}_t(dt,dx) 
 + 
\int_{\mathbb{R}}(e^{\lambda (t,T)x}-1)
\]
\[ 
(1-\frac1{\beta (t,x,T_{i(T)},T+\delta)})\nu_t^{i(T)}(dt,dx)+
 \lambda (t,T)c_t^{\frac12}dW_t^{i(t)}+\int_{\mathbb{R}}(e^{\lambda (t,T)x}-1)(\mu-\nu_t)^{i(t)})(dt,dx)  \Big{)}.\qquad \forall T\in \mathbb{R}_+. 
 \]
\end{enumerate}
That way, we get an existence result for a class of term structure models without using a priori existence of an HJM model giving rise to the continuous (and unbounded) tenor LIBOR-Market-Model.$\quad\Box$
\begin{corr}
If an HJM model exists for all $T>0$, then so does the LIBOR termstructure in theorem (\ref{LIBORTERM}) driven by the same driving process. Therefore the conditions (\ref{cond1}), (\ref{cond2}) and (\ref{cond3}), for all possible equidistant grids with arbitrary $\delta >0$, are necessary conditions on the driving process for the existence of the $HJM$ model driven by that process. 
\end{corr}
There is another interpolation method by Schloegl in \cite{SS2002} pages 197-218, which is more flexible, but does not necessarily lead to an SDE analogous to the continuously compounded risk neutrally modeled rates of \cite{BGM}.
\section{Conclusion}
\label{conclusio}
We have derived necessary conditions on the driving process for the existence of a solution to the HJM-equation. Sufficient conditions can be found in \cite{FT2008} and \cite{JZ2007}. It would be exciting to close the gap and derive real criteria for existence (calculating the actual limiting behavior for $\delta \rightarrow 0$). This may constitute further work.
\newpage

\bibliography{dissertation}

\end{document}